\documentclass[conference]{IEEEtran}
\usepackage[utf8]{inputenc}
\usepackage[english]{babel}
\usepackage{blindtext}
\usepackage{amssymb}
\usepackage{amsmath}
\usepackage{nicefrac}
\usepackage{mathtools}
\usepackage[keeplastbox]{flushend}
\usepackage{todonotes}

\DeclarePairedDelimiter\norm{\lVert}{\rVert}
\title{Using Approximate Computing for the Calculation of Inverse Matrix $p$-th Roots}
\author{
  \IEEEauthorblockN{Michael Lass}
  \IEEEauthorblockA{Department of Computer Science\\Paderborn University\\33098 Paderborn, Germany\\michael.lass@uni-paderborn.de}
  \and
  \IEEEauthorblockN{Thomas D. K\"uhne}
  \IEEEauthorblockA{Department of Chemistry\\Paderborn University\\33098 Paderborn, Germany\\tdkuehne@mail.uni-paderborn.de}
  \and
  \IEEEauthorblockN{Christian Plessl}
  \IEEEauthorblockA{Department of Computer Science\\Paderborn University\\33098 Paderborn, Germany\\christian.plessl@uni-paderborn.de}
}

\begin{document}

  \maketitle

  \begin{abstract}

    Approximate computing has shown to provide new ways to improve performance
    and power consumption of error-resilient applications. While many of these
    applications can be found in image processing, data classification or
    machine learning, we demonstrate its suitability to a problem from
    scientific computing. Utilizing the self-correcting behavior of iterative
    algorithms, we show that approximate computing can be applied to the
    calculation of inverse matrix $p$-th roots which are required in many
    applications in scientific computing. Results show great opportunities to
    reduce the computational effort and bandwidth required for the execution of
    the discussed algorithm, especially when targeting special accelerator
    hardware.

  \end{abstract}

  \section{Introduction}

  Approximate computing has gained a lot of attention over the last years as a
  technique to accelerate systems or applications or increase their power
  efficiency by allowing errors or a loss of precision. Many different
  techniques have been proposed, such as bit truncation, approximate adder
  circuits, analog computing and usage of neural networks. Typical applications
  for approximate computing are from the areas of image processing, data
  classification and machine learning where small variations in the output are
  barely noticeable and can be tolerated.

  Chippa et al.~\cite{ChippaChakradharRoyEtAl2013} showed that additionally,
  iterative algorithms are often inherently resilient to approximation because
  errors introduced in one iteration are likely to be fixed in later iterations.
  Iterative algorithms can be found for many numeric problems encountered in
  scientific computing. Klavik et al.~\cite{KlavikMalossiBekasEtAl2014}
  demonstrated that the \emph{Conjugate Gradient} method used to solve systems
  of linear equations performs well using low precision arithmetic for the
  computationally expensive parts. Schöll et al.~\cite{schoell2015,schoell2016}
  describe the application of efficient fault-tolerance to Preconditioned
  Conjugate Gradient in order to execute this algorithm on approximate hardware.
  In this work we focus on an iterative algorithm which is used to calculate the
  inverse $p$-th root $A^{-\nicefrac{1}{p}}$ for a given symmetric positive
  definite matrix $A$.

  Calculating inverse $p$-th roots is important for a large variety of different
  applications in scientific computing such as preconditioning, Kalman
  filtering, non-linear optimization, solving linear least squares problems and
  systems of linear equations,
  as well as generalized eigenvalue problems, particularly solving Schr\"odinger
  and Maxwell equations, to name just a few. While for the former applications
  an approximation is generally sufficient, the latter ought to be solved
  exactly. Nevertheless, as will be demonstrated in the following, it is even
  possible to obtain essentially exact results using the approximate computing
  paradigm.

  In this work, we contribute an analysis of the error resiliency of the given
  algorithm and discuss the potential benefits of using low-precision arithmetic
  and data representation in applications using this algorithm to calculate
  inverse $p$-th roots of large matrices.

  \section{Calculation of Inverse Matrix $p$-th Roots}

  The iterative algorithm discussed in this work was first described by Bini et
  al.~\cite{bini2005algorithms}. We briefly describe the algorithm in the
  following. Let $C_k,\enspace k=0,1,\dots$ be the sequence of intermediate
  result matrices. Starting with an initial guess $C_0$, in each iteration the
  result is refined by
  \begin{align}
    C_{k+1} = \frac{1}{p}\left( (p+1)C_k - C_k^{p+1} A \right)
  \end{align}
  If the initial guess $C_0$ was already close to $A^{-\nicefrac{1}{p}}$ such
  that
  \begin{align}
    \norm{I - C_0^p A}_2 < 1 \label{eq:constraint}
  \end{align}
  $C_k$ for $k\rightarrow\infty$ converges against $A^{-\nicefrac{1}{p}}$. For
  $p=1$ the algorithm corresponds to the well-known Newton-Schulz
  method~\cite{schulz1933} used to iteratively calculate inverse matrices.

  There are different possible choices for $C_0$. For our evaluation we use
  \begin{align}
    C_0 = {\left(\norm{A}_1 \cdot \norm{A}_\infty\right)}^{-1} A^\mathsf{T}
  \end{align}
  which is proven to always fulfil constraint~(\ref{eq:constraint}) and
  therefore guarantees convergence~\cite{richters2016}.

  \section{Use-Cases for Approximation}

  Approximate computing is applicable to the presented algorithm in different
  use-cases. As described, some applications benefit from having an
  approximation of inverse $p$-th roots if these can be calculated quickly.
  However, even if exact results are required, the iterative nature of the
  algorithm allows using imprecise arithmetic for the majority of iterations and
  then refining the result in few iterations using precise arithmetic. There are
  different scenarios that suggest using approximate arithmetic or approximate
  storage which are discussed in the following.

  \subsection{Approximate Arithmetic}
  Approximating the arithmetic operations in the presented algorithm or using
  low precision can be beneficial for overall performance or the energy
  efficiency, depending on the underlying compute platform.
  Recent GPUs targeting the data center, such as the \emph{NVIDIA Tesla P100},
  support half-precision\footnote{Half-precision: 1 sign bit, 5 exponent bits,
  10 stored mantissa bits} floating-point operations, doubling the peak
  performance compared to single-precision\footnote{Single-precision: 1 sign
  bit, 8 exponent bits, 23 stored mantissa bits} arithmetic and quadrupling it
  compared to double-precision\footnote{Double-prevision: 1 sign bit, 11
  exponent bits, 52 stored mantissa bits} arithmetic.~\cite{tesla}

  For custom hardware accelerators, the data width also plays an important role.
  Currently, FPGAs are becoming more and more common as an accelerator platform,
  next to CPUs and GPUs. Reduced precision can significantly lower resource
  requirements on FPGAs, e.g.\ for current Xilinx devices using less than 17
  stored mantissa bits reduces the number of required DSPs by half compared to
  single-precision. For custom CMOS designs it has been shown that the power
  consumption of multipliers rises at least quadratically with the number of
  input bits~\cite{callaway1997}. Because the delay increases too, this has an
  amplified effect on the energy consumption. Using fixed-point arithmetic with
  only few bits can further simplify custom designs for FPGAs or ASICs.

  Lastly, more advanced approximate computing techniques such as biased voltage
  over-scaling, approximate arithmetic circuits or specifically targeted
  overclocking have the potential to improve clock speeds or lower the
  complexity and power consumption of custom hardware designs.

  \subsection{Approximate Storage}
  In applications using the examined algorithm, data sets are typically very
  large, i.e.\ calculations need to be performed on matrices containing billions
  of entries. This not only leads to great computational demands but also
  requires large memories and high memory bandwidth. Even if calculations are
  performed precisely, the required memory space and bandwidth can be reduced by
  storing intermediate results after each iteration as well as the input data
  itself with less precision, using fewer bits.

  If the computations are to be performed on special hardware accelerators such
  as GPUs or FPGAs, the data needs to be transferred between the host and these
  devices. In this case, bandwidth quickly becomes a limiting factor, motivating
  the use of low precision representations for the transferred data.

  In the scope of this work, we evaluate the impact of both, low precision
  arithmetic and low precision storage of intermediate results, on the examined
  algorithm.

  \section{Evaluation}

  To assess the resiliency of the given algorithm to certain errors, we simulate
  simple approximation techniques. First, we use custom-precision floating-point
  and fixed-point arithmetic for \emph{all involved calculations} in order to
  get an understanding of the required data ranges and precision. Second, we
  restrict the use of these data types to the input data and stored intermediate
  results, simulating \emph{approximate storage} or \emph{data exchange} of the
  involved matrices.

  \subsection{Problem and Data Set}
\label{sec:dataset}

  We use the examined algorithm to efficiently compute two particularly
  time-consuming kernels of effective single-particle Schr\"odinger equations,
  which are to orthogonalize and solve the corresponding generalized eigenvalue
  problem. For that purpose, we first calculate the inverse square-root
  $\mathbf{S}^{-\nicefrac{1}{2}}$ of the so-called overlap matrix $\mathbf{S}$,
  which is of dimension $N \times N$ and whose elements read as $S_{ij}=\langle
  \varphi_i | \varphi_j \rangle$, where $\varphi_i$ are the $N$ non-orthogonal
  basis functions spanning the Hilbert space, and subsequently $\mathbf{S}^{-1}$
  to efficiently solve the resulting orthogonalized matrix eigenvalue
  equation~\cite{Ceriotti2008, Richters2014}.

  \begin{figure}
    \centering
    \includegraphics[width=.7\columnwidth]{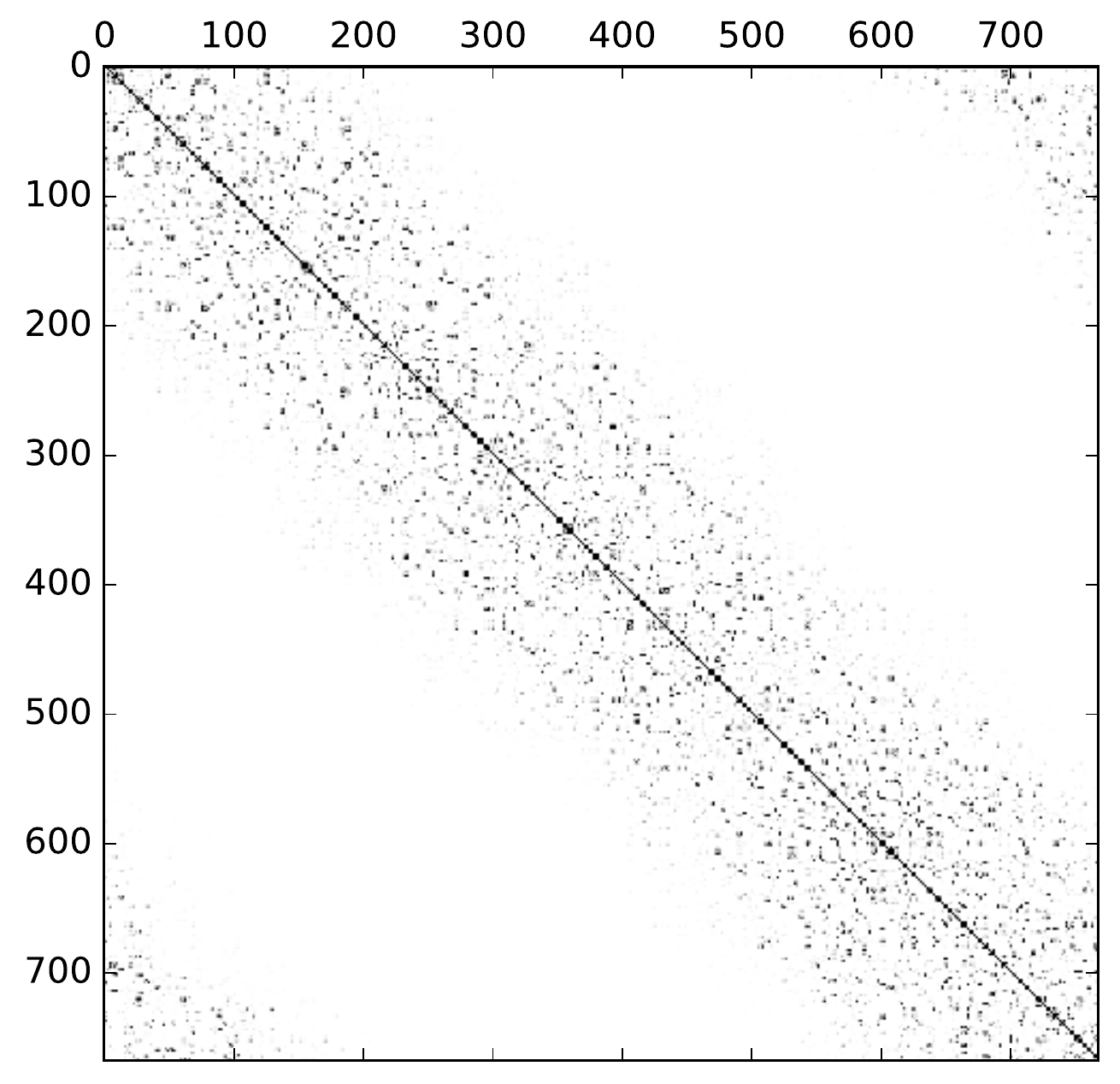}
    \caption{Structure of an examined overlap matrix $\mathbf{S}$ ($N=786$).}
\label{fig:matrix}
  \end{figure}

  Figure~\ref{fig:matrix} shows an overlap matrix $\mathbf{S}$ of size $N=786$
  which we used for our evaluation. It is a symmetric positive definite matrix
  with 25\% non-zero elements in the range of $[-1,1]$, representing a system of
  128 H\textsubscript{2}O molecules. For larger matrices density decreases to
  12.4\% ($N=1572$), 6.2\% ($N=3072$) and 3.1\% ($N=6144$) non-zero elements.
  This is a manifestation of the \emph{nearsightedness} principle of electronic
  matter and the foundation of linear scaling electronic
  algorithms.~\cite{prodan2005}

  \subsection{Methodology}

  For the presented simulations we use Python along with NumPy and
  SciPy~\cite{scipy}, which provide the required data structures and numeric
  operations. This allows us to define entirely custom data types, e.g.\
  floating-point types with a custom number of bits in the exponent and
  mantissa, and fixed-point types with a selectable number of bits and
  selectable scaling factor. Implementing basic arithmetic operations like
  addition,
  subtraction
  and multiplication
  for our custom data types enables NumPy to use these data types in its own
  array data structures and numeric operations.

  This approach allows a very flexible and fast implementation of simulators for
  different approximation techniques. Besides the mentioned floating-point and
  fixed-point data types, influences like noise or random bitflips can be easily
  implemented and adjusted. This flexibility comes at the cost of a performance
  penalty as each arithmetic operation now implies doing a function call,
  performing the necessary simulation steps and instantiating a return object.
  We deal with this performance degradation by implementing all classes in
  Cython~\cite{cython}, producing statically typed C-code which executes orders
  of magnitude faster than interpreted Python code, and distributing the
  simulations over many machines of a large compute cluster.

  \subsection{Results}

  \subsubsection{Overall Error Resiliency} To assess the overall error
  resiliency of the algorithm, we chose overlap matrices of dimension $N=768$
  and set $p=2$ to calculate the inverse square root for these matrices, which
  is one step of solving the generalized eigenvalue problem as described in
  Section~\ref{sec:dataset}. With simulation, we determine the convergence of
  the algorithm, depending on the given precision. The iterative algorithm shows
  to be rather resilient to low precision, both for storage of intermediate
  matrices as well as for all used arithmetic operations.

  \begin{figure}
    \centering
    \includegraphics[width=\columnwidth]{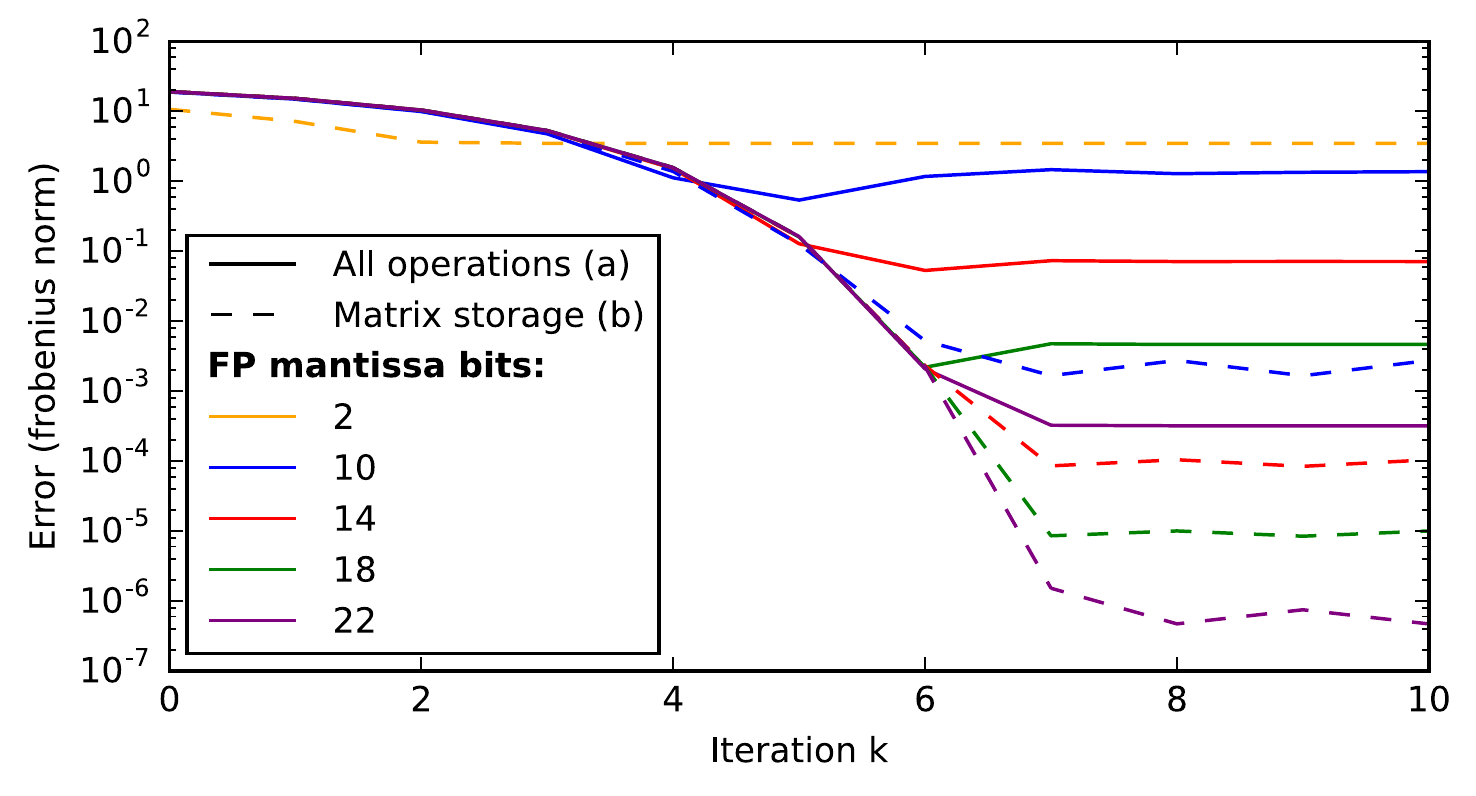}
    \caption{Convergence behavior when using custom-precision floating-point for (a) all arithmetic operations and (b) only for storage of intermediate results.}
\label{fig:arbFloat}
  \end{figure}

  Figure~\ref{fig:arbFloat} shows the error between the intermediate solutions
  obtained from the algorithm using floating-point with custom mantissa widths
  and a solution that was precomputed using double-precision. The error is
  determined by the Frobenius norm

  \begin{align*}
    \norm{C_k - \mathbf{S}^{-\nicefrac{1}{p}}}_F \coloneqq \sqrt{\sum_{i=1}^N \sum_{j=1}^N \vert \gamma_{ij} - \chi_{ij}\vert^2}
  \end{align*}
  where $C_k = (\gamma_{ij})$ and $\mathbf{S}^{-\nicefrac{1}{p}} = (\chi_{ij})$.

  The observed convergence of the algorithm can be split into two phases: First
  the error steadily decreases, following a curve described by
  \begin{align*}
    f(x) = 2^{-P(x)}
  \end{align*}
  with $P(x)$ being a polynomial function of at least second degree. In the
  second phase, being limited by the given precision of the data type, the
  algorithm does not converge further but oscillations may be observed.

  This shows that the convergence in the first phase is barely influenced by the
  introduced errors. Only for less than 10 mantissa bits the algorithm does not
  converge at all. Consequently, half-precision floating-point arithmetic is
  sufficient to retain convergence. Approximation does however increase the
  lower bound for the error. Therefore the second phase of conversion start
  earlier for lower precision. Increasing the precision in later iterations
  allows the algorithm to further converge against a lower error, opening the
  possibility for dynamic precision scaling. Observing the changes introduced in
  each iteration, the necessity of increased precision can be detected at
  runtime.

  Approximating only the storage of intermediate results allows significantly
  stronger approximation while achieving similar precision in the output. For
  example, storing only 10 mantissa bits allows a similar error as doing all
  calculation using 18 mantissa bits. In our use case the algorithm still
  converges if only two mantissa bits are used for all stored values.

  \begin{figure}
    \centering
    \includegraphics[width=\columnwidth]{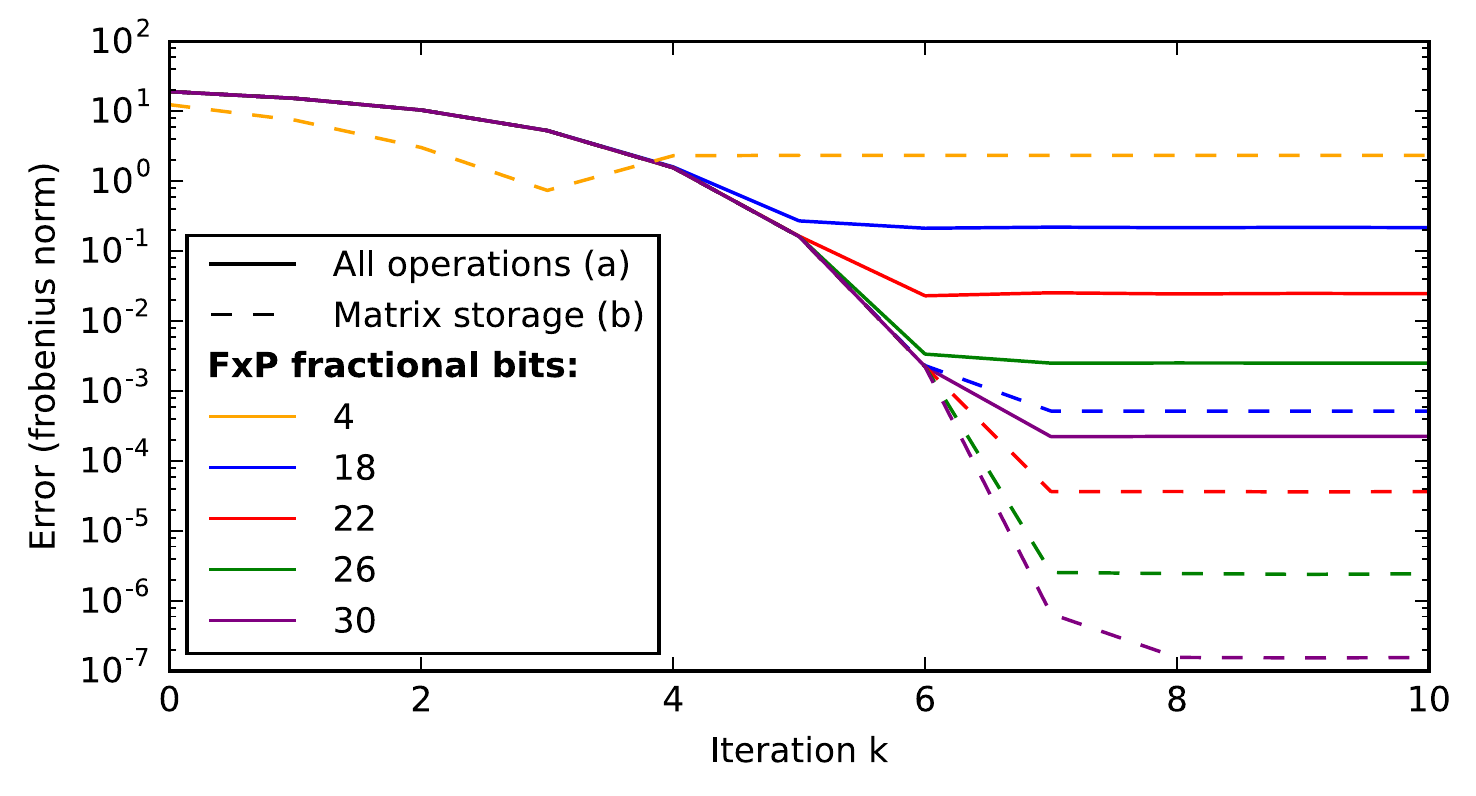}
    \caption{Convergence behavior when using custom-precision fixed-point for (a) all arithmetic operations and (b) only for storage of intermediate results.}
\label{fig:fixedpoint}
  \end{figure}

  Figure~\ref{fig:fixedpoint} shows similar behavior when using fixed-point
  arithmetic with a low number of bits. To retain convergence, 18 fractional
  bits are required for arithmetic operations. Again, restricting the
  approximation to stored intermediate results permits stronger approximation.
  In our evaluation, storing only four fractional bits showed to be sufficient
  to retain convergence.

  \subsubsection{Influence of the Matrix Size}

  \begin{figure}
    \centering
    \includegraphics[width=\columnwidth]{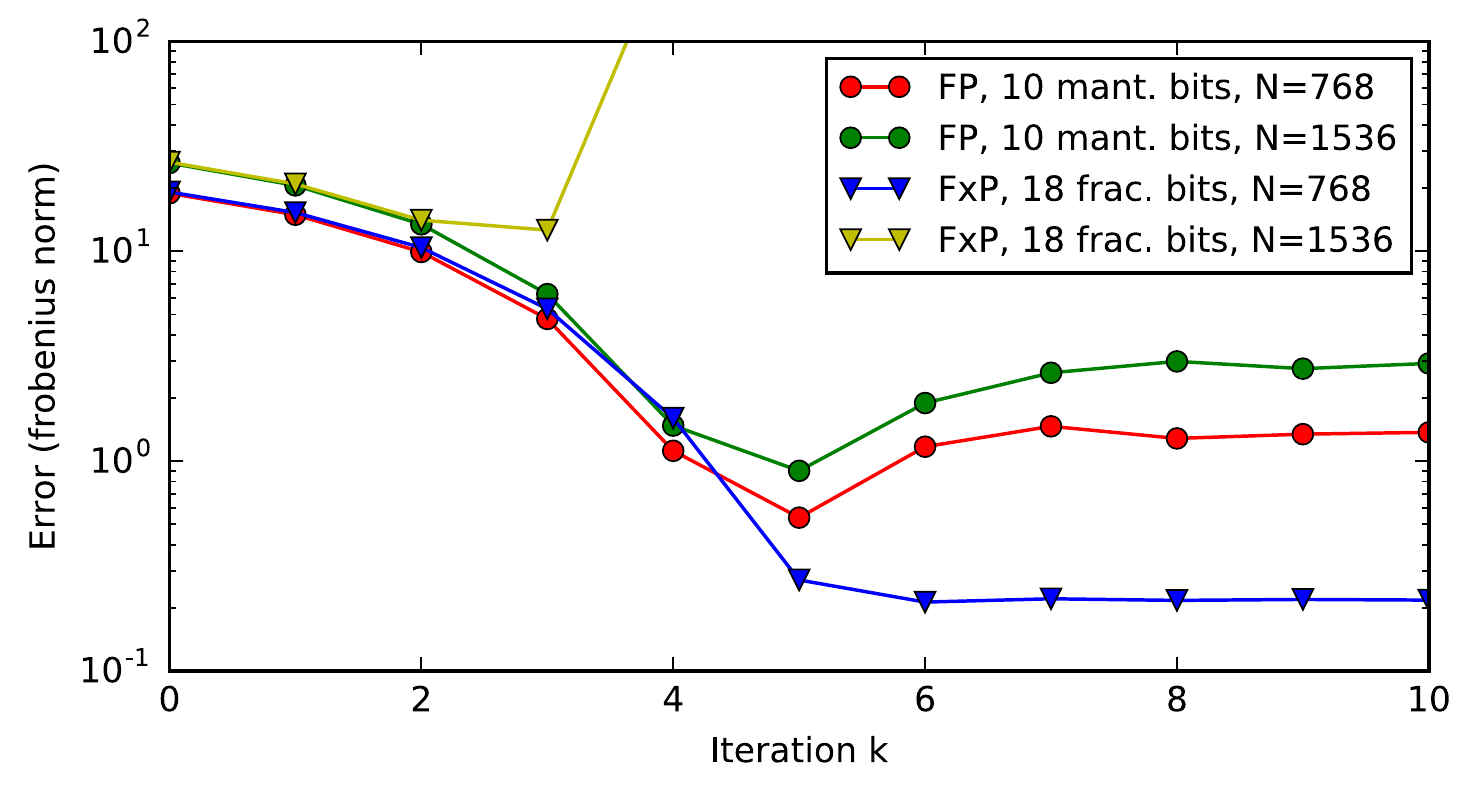}
    \caption{Convergence behavior using custom-precision floating-point (FP) and fixed-point (FxP) for different matrix sizes.}
\label{fig:sizes}
  \end{figure}

  Most of the results presented before apply directly to larger matrices from
  our problem set, in particular when only approximating the storage of
  intermediate results. Approximating all arithmetic operations using
  low-precision fixed-point arithmetic however exhibits a limitation. As shown
  in Figure~\ref{fig:sizes}, using 18 fractional bits is sufficient to retain
  convergence for $N=786$ but for $N=1536$ the error eventually increases. The
  reason for this behavior is that larger matrices from our set are more sparse
  (see Section~\ref{sec:dataset}) and therefore their inverse contain smaller
  values which cannot be represented appropriately with the given number of
  fractional bits.

  For floating-point this effect is not relevant, as depicted in
  Figure~\ref{fig:sizes}. The slightly larger final error for $N=1536$ can be
  explained by the use of the Frobenius norm as error metric since it adds up
  the quadratic errors of all matrix elements.

  \subsubsection{Influence of $p$} 
\begin{figure}
    \centering
    \includegraphics[width=\columnwidth]{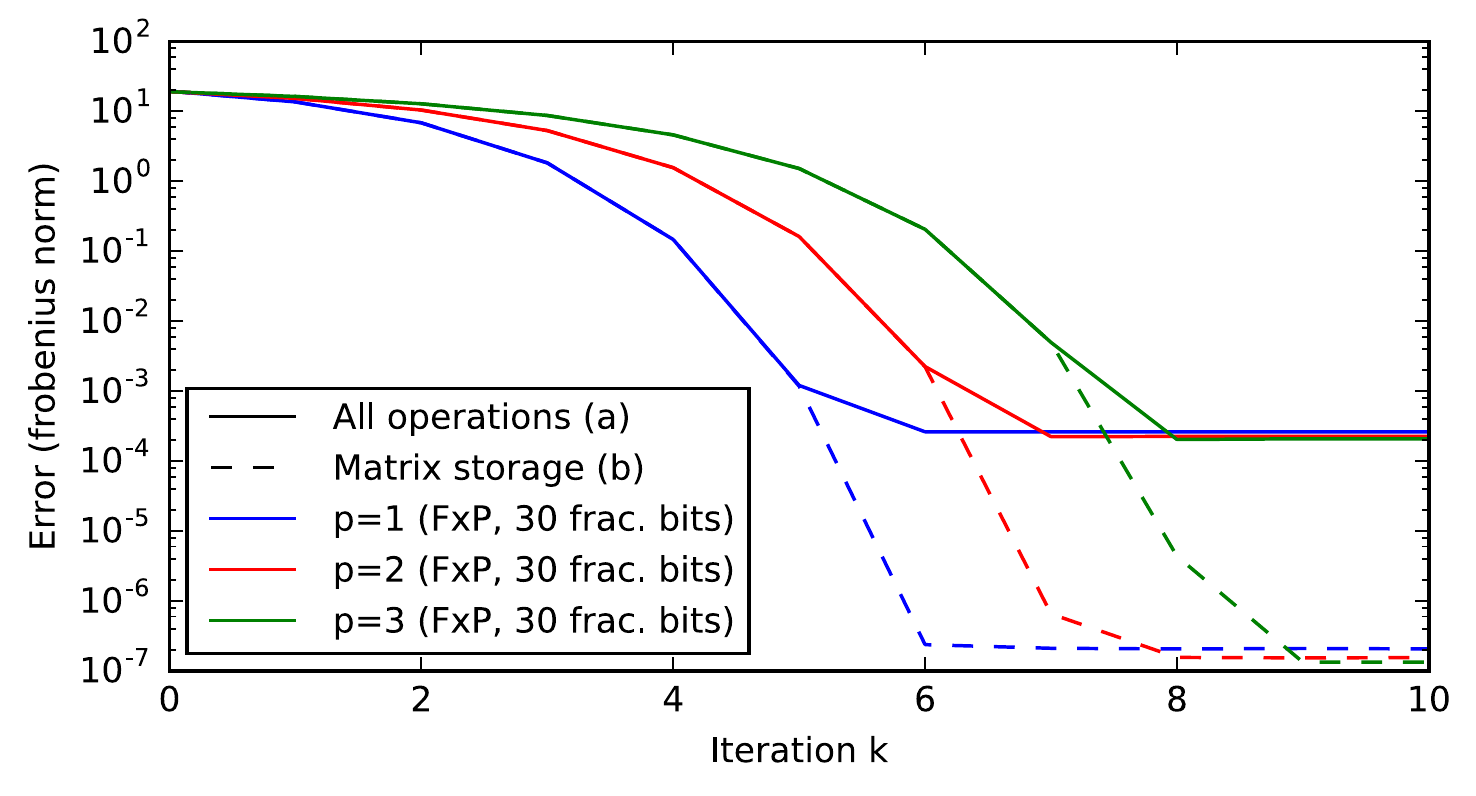}
    \caption{Convergence behavior using custom-precision fixed-point (FxP) for different $p$.}
\label{fig:diff-p}
  \end{figure}

Calculating the inverse $p$-th root for
  $p\neq2$ shows similar behavior as for $p=2$, as shown in
  Figure~\ref{fig:diff-p} for custom-precision fixed-point arithmetic and
  storage. With increasing $p$ the algorithm in general needs an increasing
  number of iterations to converge. This effect is independend of the applied
  approximation.

  \section{Conclusion}

  The presented results show the resiliency of the examined algorithm against
  errors introduced due to low precision arithmetic and storage. While a certain
  precision has to be provided to retain convergence of the algorithm, further
  precision is only required in final iterations if a precise result is desired.

  It stands out that the number of iterations required to reach a certain
  precision does not significantly increase with the amount of approximation.
  This sets the examined algorithm apart from other iterative methods like the
  preconditioned conjugate gradient method which was modified to run on
  approximate hardware by Schöll et al.~\cite{schoell2016} and showed to require
  additional iterations when using approximation.

  This opens up great opportunities for the acceleration of applications in the
  scientific computing domain requiring the calculation of inverse $p$-th roots:
  Using half-precision floating-point in the first iterations can lead to a
  $2\times$ speedup for these iterations on suitable GPUs. Resource requirements
  on FPGAs can be reduced by half and for custom CMOS designs, power consumption
  of the multipliers can be reduced by $\nicefrac{3}{4}$. Moreover, the overhead
  required for data exchange when using GPUs or custom hardware can be
  significantly reduced as data can be represented using low precision data
  types.

  If a precise solution for $A^{-\nicefrac{1}{p}}$ is required for the
  application, results obtained using approximation can be refined into a
  precise solution in very few additional iterations. In a scenario using
  approximate hardware accelerators, this can be done in software while leaving
  the main part of the work to the external accelerators.

  \bibliographystyle{IEEEtran}
  \bibliography{literature}

\end{document}